\magnification=1200

{\mathcode`A="7041 \mathcode`B="7042 \mathcode`C="7043 \mathcode`D="7044
\mathcode`E="7045 \mathcode`F="7046 \mathcode`G="7047 \mathcode`H="7048
\mathcode`I="7049 \mathcode`J="704A \mathcode`K="704B \mathcode`L="704C
\mathcode`M="704D \mathcode`N="704E \mathcode`O="704F \mathcode`P="7050
\mathcode`Q="7051 \mathcode`R="7052 \mathcode`S="7053 \mathcode`T="7054
\mathcode`U="7055 \mathcode`V="7056 \mathcode`W="7057 \mathcode`X="7058
\mathcode`Y="7059 \mathcode`Z="705A

\font\tencyr=wncyr10
\font\sevencyr=wncyr7
\font\fivecyr=wncyr5
\newfam\cyrfam
\def\cyr{\fam\cyrfam\tencyr}
\textfont\cyrfam=\tencyr
\scriptfont\cyrfam=\sevencyr
\scriptscriptfont\cyrfam=\fivecyr

\def\Sha{{\cyr W}}

\def\Z{{\bf Z}}
\def\Q{{\bf Q}}

\def\C{{\bf C}}
\def\P{{\bf P}}
\def\F{{\bf F}}

\def\lim{\mathop{\rm lim}\displaylimits}
\def\ord{\mathop{\rm ord}\nolimits}
\def\Card{\mathop{\rm Card}\nolimits}
\def\rk{\mathop{\rm rk}\nolimits}
\def\mod{\mathop{\rm mod.}\nolimits}

 \font\eightrm=cmr8
\def\pc#1{\rm #1\eightrm}

\newcount\refno
\long\def\ref#1:#2<#3>{                                        
\global\advance\refno by1\par\noindent                              
\llap{[{\bf\number\refno}]\ }{#1} \pointir{\it #2} #3\medbreak }

\def\citer#1{[{\bf\number#1}]}
\def\pointir{\unskip . --- \ignorespaces}

\newbox\bibbox
\setbox\bibbox\vbox{\bigbreak
\centerline{{\pc  BIBLIOGRAPHIC REFERENCES}}
\medskip

\ref{\pc BIRCH} (Brian) and {\pc SWINNERTON-\pc DYER} (Peter):
Notes on elliptic curves\/ {\rm II},
<J. Reine Angew.\ Math.\ {\bf 218} (1965), 79--108.>
\newcount\bsd \global\bsd=\refno

\ref{\pc CAPORASO} (Lucia), {\pc HARRIS} (Joe) and {\pc MAZUR} (Barry):
Uniformity of rational points,
<J. Amer.\ Math.\ Soc.\ {\bf 10} (1997) 1, 1--35.>
\newcount\cahama \global\cahama=\refno

\ref{\pc CORNELL} (Gary) and {\pc SILVERMAN} (Joseph) (Eds.):
Arithmetic geometry,
<Springer-Verlag, 1986.>
\newcount\faltings \global\faltings=\refno

\ref{\pc CORNELL} (Gary), {\pc SILVERMAN} (Joseph) and {\pc STEVENS}
(Glenn) (Eds.):
Modular forms and Fermat's last theorem,
<Springer-Verlag, 1997.>
\newcount\flt \global\flt=\refno

\ref{\pc EDIXHOVEN} (Bas):
Rational elliptic curves are modular (after Breuil, Conrad, Diamond
and Taylor),  
<Ast\'erisque {\bf 276} (2002), 161--188.>
\newcount\basmodular \global\basmodular=\refno

\ref{\pc EDIXHOVEN} (Bas):
Rational torsion points on elliptic curves over number fields
            (after Kamienny and Mazur),
<Ast\'erisque, {\bf 227} (1995) 4, 209--227.>
\newcount\basmerel \global\basmerel=\refno

\ref{\pc {}VAN \pc FRANKENHUYSEN} (Machiel):
The $ABC$ conjecture implies Vojta's height inequality for curves,
<J. Number Theory {\bf 95} (2002) 2, 289--302.>
\newcount\abcvojta \global\abcvojta=\refno

\ref{\pc KATO} (Kazuya) and {\pc TRIHAN} (Fabien):
Conjectures of Birch and Swinnerton-Dyer in positive characteristics
assuming the finiteness of the Tate-Shafarevich group,
<Preprint.>
\newcount\katotrihan \global\katotrihan=\refno

\ref{\pc KOBLITZ} (Neal):
Introduction to elliptic curves and modular forms,
<Springer-Verlag, 1993.>
\newcount\koblitz \global\koblitz=\refno

\ref{\pc LANG} (Serge):
Number theory\/ {\rm III}, Diophantine geometry, 
<Springer-Verlag, 1991.>
\newcount\lang \global\lang=\refno

\ref{\pc MANIN} (Yuri):
Cyclotomic fields and modular curves,
<Uspehi Mat.\ Nauk, {\bf 26} (1971) 6, 7--71.>
\newcount\manin \global\manin=\refno

\ref{\pc MAZUR} (Barry):
Arithmetic on curves,
<Bull.\ Amer.\ Math.\ Soc.\ (N.S.) {\bf 14} (1986) 2, 207--259.>
\newcount\mazurarithm \global\mazurarithm=\refno

\ref{\pc MAZUR} (Barry):
On the passage from local to global in number theory,
<Bull.\ Amer.\ Math.\ Soc.\ (N.S.) {\bf 29} (1993) 1, 14--50.>
\newcount\mazurlocal \global\mazurlocal=\refno

\ref{\pc MAZUR} (Barry):
Questions about powers of numbers,
<Notices of the Amer.\ Math.\ Soc., February 2000.>
\newcount\mazurabc \global\mazurabc=\refno

\ref{\pc MEREL} (Lo{\"{\i}}c):
Bornes pour la torsion des courbes elliptiques sur les corps de
nombres, 
<Invent.\ Math.\ {\bf 124} (1996) 1-3, 437--449.>
\newcount\merel \global\merel=\refno

\ref{\pc MUMFORD} (David):
Abelian varieties,
<Oxford University Press, 1970.>
\newcount\mumford \global\mumford=\refno

\ref{\pc NEKOV{\'A}{\v R}} (Jan):
On the parity of ranks of Selmer groups\/ {\rm II},
<C. R. Acad.\ Sci.\ Paris S\'er. I Math.\ {\bf 332} (2001) 2, {99--104}.>
\newcount\nekovar \global\nekovar=\refno

\ref{\pc N{\'E}RON} (Andr{\'e}):
Probl\`emes arithm\'etiques et g\'eom\'etriques rattach\'es \`a
la notion de rang d'une courbe alg\'ebrique dans un corps,
<Bull.\ Soc.\ Math.\ France {\bf 80} (1952), 101--166.>
\newcount\neron \global\neron=\refno

\ref{\pc OESTERL{\'E}} (Joseph):
Nouvelles approches du ``th\'eor\`eme'' de Fermat,
<Ast\'erisque {\bf 161-162} (1989), 165--186.>
\newcount\josephfermat \global\josephfermat=\refno

\ref{\pc OESTERL{\'E}} (Joseph):
Empilements de sph\`eres,
<Ast\'erisque {\bf 189-190} (1990), 375--397.>
\newcount\josephspheres \global\josephspheres=\refno

\ref{\pc PERRIN-\pc RIOU} (Bernadette):
Travaux de Kolyvagin et Rubin,
<Ast\'erisque {\bf 189-190} (1990), 69--106.>
\newcount\perrinriou \global\perrinriou=\refno

\ref{\pc RUBIN} (Karl) and {\pc SILVERBERG} (Alice):
Ranks of elliptic curves,
<Bull.\ Amer.\ Math.\ Soc.\ (N.S.) {\bf 39} (2002) 4, 455--474.>
\newcount\rubinrank \global\rubinrank=\refno

\ref{\pc SAMUEL} (Pierre):
Lectures on old and new results on algebraic curves,
<Tata Institute of Fundamental Research, 1966.>
\newcount\samuel \global\samuel=\refno

\ref{\pc SERRE} (Jean-Pierre):
A course in arithmetic,
<Springer-Verlag, 1973.>
\newcount\serrearithm \global\serrearithm=\refno

\ref{\pc SERRE} (Jean-Pierre):
Galois cohomology,
<Springer-Verlag, 2002.>
\newcount\serregalois \global\serregalois=\refno

\ref{\pc SERRE} (Jean-Pierre):
Points rationnels des courbes modulaires $X_{0}(N)$ [d'apr\`es
            Barry Mazur],
<Lecture Notes in Math. {\bf 710}, 89--100>
\newcount\serremazur \global\serremazur=\refno

\ref{\pc SERRE} (Jean-Pierre):
Lectures on the Mordell-Weil theorem,
<Friedr.\ Vieweg \& Sohn, 1997.>
\newcount\lectures \global\lectures=\refno

\ref{\pc SILVERMAN} (Joseph):
The arithmetic of elliptic curves,
<Springer-Verlag, 1986.>
\newcount\silverman \global\silverman=\refno

\ref{\pc SILVERMAN} (Joseph) and {\pc TATE} (John):
Rational points on elliptic curves,
<Springer-Verlag, 1992.>
\newcount\silvermantate \global\silvermantate=\refno

\ref{\pc STEIN} (William):
There are genus one curves over $\Q$ of every odd index,
<J.\ Reine Angew.\ Math.\ {\bf 547} (2002), 139--147.>
\newcount\stein \global\stein=\refno

\ref{\pc TATE} (John):
On the conjectures of Birch and Swinnerton-Dyer and a geometric
analog\/{\rm [ue]},
<S\'eminaire Bourbaki, Vol.\ 9, Exp.\ 306, Soc. Math. France, 1995.>
\newcount\tatebsd \global\tatebsd=\refno

\ref{\pc ULMER} (Douglas):
Elliptic curves with large rank over function fields,
<Ann.\ of Math.\ (2) {\bf 155} (2002) 1, 295--315.>
\newcount\ulmer \global\ulmer=\refno

\ref{\pc WEIL} (Andr{\'e}):
Number theory, an approach through history,
<Birkh\"auser, 1984.>
\newcount\weilnt \global\weilnt=\refno

}

{\openup,5\jot
\advance\parskip by3pt plus2pt

\hbox{\vbox{\kern2cm}}
\centerline{\bf Arithmetic on curves}

\smallskip

\centerline{Chandan Singh {\pc DALAWAT}}

\bigskip

We intend to give a brief account of what is known or conjectured
about the set $C(\Q)$ of rational points on a smooth projective
absolutely connected curve $C$ of genus $g$ over $\Q$.  The idea is to
show how the arithmetic properties of algebraic curves are governed by
the familiar trichotomy~: $g=0$, $g=1$, $g\ge2$.  Only incidentally
shall we mention fields other than $\Q$ and varieties other than
curves.

We shall not give even the definitions of all the concepts which
appear below, let alone the proofs of the theorems we state.  Our wish is
merely to provide some references --- for the most part to survey
articles and books --- which the interested reader is advised to look
up for the precise definitions, statements and --- where they are
known --- proofs.

If I were to summarise in a few lines the work of more than a hundred
mathematicians on this topic over the last century, I'd say~: if
$g=0$, then $C(\Q)$ is either empty or infinite~; we know how to
distinguish between the two cases and, in the latter case, how to
describe it completely.  If $g=1$, we don't know when $C(\Q)$ can be
empty~; if it is not empty, it is a torsor under a finitely generated
commutative group~$T$ attached to $C$~; we know how to compute the
torsion subgroup of $T$ but not its rank~; we know which groups can
occur as the torsion but not which integers as ranks.  Finally, if
$g\ge2$, we know that $C(\Q)$ is finite but not how ``large'' its
points can be, nor whether the number of points remains bounded as $C$
varies among curves of genus~$g$.  These cryptic remarks will now be
slightly amplified.

A smooth projective curve of genus~0 can be realised as a plane conic
$C$ of the form $aX^2+bY^2+cZ^2=0$ ($a,b,c\in\Q^\times$).  A theorem
of Legendre \citer\weilnt\ can be rephrased as saying that for $C$ to
have a rational point it is sufficient for it to have a real point ---
which amounts to saying that $a,b,c$ are not of the same sign --- and,
for every odd prime $p$, a point over $\Q_p$.  As a matter of fact, it
is sufficient to check solvability at the real place and --- after
rescaling $X,Y,Z$ so that the coefficients $a,b,c$ are in $\Z$ --- to
check that the congruence $aX^2+bY^2+cZ^2\equiv0\mod M$ has a
primitive solution for a certain integer $M$ determined by the
coefficients.  Thus, checking for solvability is a finite amount of
computation~; so is finding a solution if there is one.  If
$C(\Q)=\hbox{\O}$, one can write down a quadratic extension~$K$ of
$\Q$ such that $C(K)\neq\hbox{\O}$.  If $C$ has a rational point, then
it is isomorphic to the projective line $\P_{\!1}$, so that we have an
explicit description of the set $C(\Q)$ of all its rational points.
Similar results continue to hold over finite extensions of $\Q$.

The set $\P_{\!1}(\Q)$ comes equipped with a natural ``\thinspace
height\thinspace'' function $H$ which measures the size of a rational
point $P$~: one defines $H(P)=\sup_{i}|x_i|$ for $P=(x_0:x_1)$, where
$x_0$, $x_1$ are rational integers with no common factor (determined
by $P$ up to sign).  The number of points of height $\le B$ grows as
$\displaystyle {2\over\zeta(2)}B^2$ when $B\rightarrow+\infty$.  There
is also a good expression for the ``\thinspace error term\thinspace''
\citer\lectures.

Curves of genus~0 are the 1-dimensional case of quadrics --- smooth
projective varieties defined by a quadratic form.  Hasse's theorem
asserts that the local to global principle continues to hold for
quadrics~: if a smooth quadric has a point over every completion of
$\Q$, then it has a rational point \citer\serrearithm.  Curves of
genus~0 are also the 1-dimensional case of varieties which are {\it
potentially\/} isomorphic to the projective space $\P_{\!n}$~: those
varieties which become isomorphic to $\P_{\!n}$ over
$\overline\Q$ --- one also says that they are
twisted forms of $\P_{\!n}$.  The local to global principle holds for
them as well.  These results about smooth quadrics and about twisted
forms of projective spaces remain valid over any number field~$K$
\citer\serregalois. 

The asymptotic growth of the number of points of bounded height on
$\P_{\!n}(K)$ has been established by Schanuel, although the best
possible ``\thinspace error term\thinspace'' seems to be known only in
the case of the projective line over $\Q$ \citer\lectures.

\medbreak

In contrast to the genus~0 case, the local to global principle fails
for curves of genus~1.  The best-known example is the curve
$3X^3+4Y^3+5Z^3=0$ which has points in every completion of $\Q$ but no
rational points \citer\mazurlocal.  There is no proven algorithm to
decide, given a curve $C$ of genus~1, whether $C(\Q)$ is empty or not.
If $C(\Q)=\hbox{\O}$, there appears to be no characterisation of the
finite extensions~$K$ of $\Q$ such that $C(K)\neq\hbox{\O}$ or of the
least integer $n\ge0$ which can be the degree of a $0$-cycle on~$C$
\citer\stein.

If $C(\Q)$ is not empty and if we fix a point $O\in C(\Q)$, there is a
unique (commutative) group law $C\times C\rightarrow C$ for which $O$
is the neutral element~; a curve with such a group law is called an
{\it elliptic curve\/} \citer\silverman\ --- it would have been more
appropriate to call them {\it abelian curves\/} since they are
1-dimensional abelian varieties.  A theorem of Mordell asserts that
the commutative group $E(\Q)$ is {\it finitely generated\/}
(\citer\silverman, \citer\silvermantate) for any elliptic curve $E$
over $\Q$.  There is no proven algorithm to decide, given an elliptic
curve $E$ over $\Q$, whether the group $E(\Q)$ is finite or not.  If
$E(\Q)$ is finite, there seems to be no characterisation of those
finite extensions~$K$ of $\Q$ for which $E(K)$ is infinite.

For a given elliptic curve $E$ over $\Q$, a fairly easy theorem of
{\'E}lisabeth Lutz and Nagell determines the torsion subgroup of
$E(\Q)$ \citer\silverman.  A deep theorem of Mazur
\citer\serremazur\ gives the list of all the (finite commutative)
groups which can occur as the torsion subgroup of an elliptic curve
over $\Q$~; this list is finite and explicitly given~; given a group
$G$ on the list, one knows the elliptic curves of which $G$ is the
rational torsion.  On the other hand, given an elliptic curve $E$ over
$\Q$, it is not easy to compute the rank of the group $E(\Q)$ and
there is no proven algorithm which is guaranteed to do so.  It is not
known which numbers can occur as the rank of $E(\Q)$ when $E$ varies
among all elliptic curves defined over $\Q$
\citer\rubinrank.

In any attempt at computing the group $E(\Q)$, one comes across the
group $\Sha(E,\Q)$ of $E$-torsors over $\Q$ which are trivial over
every completion of $\Q$~; it is widely conjectured to be finite.
Manin has proved that its finiteness for every elliptic curve would
allow us to give an algorithm to test the existence of a rational
point on a curve of genus~1 \citer\lang.

The most important open problem in the arithmetic theory of
elliptic curves is the conjecture of Birch and Swinnerton-Dyer
\citer\bsd.  It involves the $L$-function $L(E,s)$ of $E$  --- a
function of a complex variable~$s$ defined using the number of points
of $E$ over the various finite fields $\F_p$ (this can be made
precise) \citer\silverman.  It was widely conjectured that $L(E,s)$
admits an analytic continuation to the whole of $\C$.  This is now
known as a consequence of the seminal work of Wiles, continued by
Breuil, Conrad, Diamond and Taylor
\citer\flt, \citer\basmodular.  Birch and Swinnerton-Dyer
conjecture that {\it the rank of\/ $E(\Q)$ equals the order of
vanishing\/ $r=\ord_{s=1}L(E,s)$} of $L(E,s)$ at $s=1$.  Assuming the
finiteness of the group $\Sha(E,\Q)$, a refined version of this
conjecture gives an expression for the ``special value''
$\displaystyle\lim_{s\rightarrow1} {L(E,s)^{\phantom{r}}\over
(s-1)^r}$ in terms of $\Card\Sha(E,\Q)$ and certain other arithmetic
invariants of~$E$.  Manin has proved that the truth of this conjecture
would give an effective method for computing the order of $\Sha(E,\Q)$
and a system of generators for $E(\Q)$ \citer\manin.

\medbreak

Results of Gross \& Zagier and Kolyvagin, among others, imply that
{\it if\/ $L(E,s)$ has at most a simple zero at\/ $s=1$, then\/
$\Sha(E,\Q)$ is finite,\/ $\rk E(\Q)=\ord_{s=1} L(E,s)$ and the
conjectured formula for the special value of\/ $L(E,s)$ at\/ $s=1$ is
true} --- to within small powers of small primes \citer\mazurlocal,
\citer\perrinriou.  Practically nothing is known about the rank of
$E(\Q)$ when the order of vanishing of $L(E,s)$ is $\ge2$ (in fact, no
elliptic curve is known whose $L$-function can be shown to vanish to
order $\ge4$).  In the other direction, it is not known that if
$E(\Q)$ is finite, then $L(E,1)\neq0$ or that if $E(\Q)$ is of rank~1,
then $L(E,s)$ has a simple zero at $s=1$.  An important recent result
of Nekov{\'a}{\v{r}} \citer\nekovar\ says that {\it if\/ $\Sha(E,\Q)$
is finite, then\/ $\ord_{s=1} L(E,s)$ and\/ $\rk E(\Q)$ have the same
parity}.

Which integers $d\ge1$ can be the area of a right triangle with
rational sides~?  This ancient problem comes down to determining the
$d$ for which the group $E_d(\Q)$ is infinite, $E_d$ being the
elliptic curve $dy^2=x^3-x$.  Tunnell has given a characterisation
under the assumption that $E_d(\Q)$ is infinite if $L(E_d,1)=0$ --- a
very special case of the Birch and Swinnerton-Dyer conjecture.  Even
this particular case is not yet known \citer\koblitz.

Mazur's theorem about the possible torsion subgroups has been extended
by Merel \citer\merel\ to all number fields, although the explicit
list of all the possibilities is known only in a few cases
\citer\basmerel.

\medbreak

Elliptic curves are the 1-dimensional case of abelian varieties.  The
theorem of Mordell has been generalised to abelian varieties over
number fields by Weil \citer\mumford\ and indeed to abelian varieties
over fields finitely generated over their prime subfields by N{\'e}ron
\citer\neron.  He has also given an expression for the asymptotic
growth of the number of points of bouded height on a projectively
embedded abelian variety over a number field \citer\lectures.

There is a version \citer\tatebsd\ of the conjecture of Birch and
Swinnerton-Dyer for abelian varieties $A$ over global fields $K$ which
predicts --- assuming the analytic continuation of the $L$-function
$L(A,s)$ to the whole of $\C$ --- that the rank of the (finitely
generated commutative) group $A(K)$ of $K$-rational points of $A$ is
equal to $r=\ord_{s=1} L(A,s)$.  There is also a conjectural
expression for the ``special value''
$\displaystyle\lim_{s\rightarrow1} {L(A,s)^{\phantom{r}}\over
(s-1)^r}$ in terms of certain arithmetic invariants of $A$ over $K$,
among them the order of the conjecturally finite group $\Sha(A,K)$ of
$A$-torsors over~$K$ which are trivial over every completion of $K$.

Elliptic curves $E$ defined over the function field $K$ of a curve
over a finite field $k$ have been used by Elkies and Shioda to give
the best known examples of sphere packings in certain dimensions
\citer\josephspheres.  Ulmer \citer\ulmer\ has given examples of
$E\,|\,K$ which have arbitrarily high ranks and such that
$E_{\overline K}$ is not definable over $\overline k$.  In both these
results, the known cases of the Birch and Swinnerton-Dyer conjecture
play a prominent role.

An important recent result of Kato and Trihan \citer\katotrihan\ says
that {\it if the group\/ $\Sha(A,K)$ is finite for an abelian
variety\/ $A$ over a function field\/ $K$ over a finite field, then
the conjecture of Birch and Swinnerton-Dyer is true for $A$ }: the
rank of $A(K)$ equals the order $r$ of vanishing of $L(A,s)$ at $s=1$
and equality holds in the conjectured expression for
$\displaystyle\lim_{s\rightarrow1} {L(A,s)^{\phantom{r}}\over
(s-1)^r}$.  

Let us now come to a curve $C$ of genus $\ge2$ over a number field
$K$.  Shuji Saito, extending Manin's result for the case $g=1$, has
proved that the failure of the local to global principle for the
existence on $C$ of a $0$-cycle of degree~1 can be accounted for by an
obstruction introduced by Manin if the group $\Sha(J,K)$ --- where $J$
is the jacobian of $C$ --- is finite
\citer\lang.  

\medbreak

The most striking result about the set $C(K)$ is Faltings' proof
\citer\faltings\ that it is {\it finite}, as had been conjectured by
Mordell~; the function field case was treated earlier by Grauert,
Manin and Samuel \citer\samuel.  Let $f\in\Q[X,Y]$ be an absolutely
irreducible polynomial defining a (smooth projective) curve $C$ of
genus $\ge2$.  An important open problem asks for a bound on the size
or height a point in $C(\Q)$ in terms of the coefficients of $f$.  It
has been proved \citer\abcvojta\ that a version of the celebrated
$abc$ conjecture (\citer\mazurabc, \citer\josephfermat) implies such a
bound.

Lang has conjectured \citer\lang\ that for a (smooth projective)
variety $V$ of general type over a number field $K$, the set $V(K)$ is
never dense in $V$~; in the case $\dim V=1$, this is a restatement of
Mordell's conjecture as proved by Faltings.  This conjecture has been
shown to imply the boundedness of $\Card C(K)$ as $C$ runs through
smooth projective curves of a given genus~$g\ge2$ over the number field
$K$ \citer\cahama.
\medbreak
Je tiens {\`a} remercier Fabien Trihan pour m'avoir fait parvenir un
exemplaire de \citer\katotrihan.
\bigbreak
}
\unvbox\bibbox
}
\bye